\documentclass[11pt]{article}
\usepackage{}
\usepackage[total={6in, 9in}]{geometry}
 \usepackage{amsfonts}
\usepackage{amsthm}
\usepackage{amssymb}
\usepackage{amsmath,amsfonts}
\usepackage{mathrsfs}
\usepackage{cases}
\usepackage{latexsym,bm}
\usepackage{indentfirst}
\usepackage{color}
\usepackage{ifpdf}
\usepackage{graphicx}
\usepackage{psfrag}
\usepackage{enumerate}
\usepackage{dsfont}
\usepackage[
pdfauthor={},
pdftitle={},
pdfstartview=XYZ,
bookmarks=true,
colorlinks=true,
linkcolor=blue,
urlcolor=blue,
citecolor=blue,
bookmarks=true,
linktocpage=true,
hyperindex=true
]{hyperref}

\usepackage[natural]{xcolor}
\usepackage{tikz}
\usepackage{subfigure,tikz}
\usepackage{tikz-3dplot}
\usepackage{pgf,tikz,pgfplots}
\usepackage{mathrsfs}

\newtheorem{THM}{\textbf{Theorem}}

\newtheorem{CLA}{\textbf{Claim}}[section]
\newtheorem{CON}[THM]{\textbf{Conjecture}}

\newcommand{\arxiv}[1]{\href{http://arxiv.org/abs/#1}{\texttt{arXiv:#1}}}

\newcommand{\qqed}{\hfill $\blacksquare$\vspace{1mm}}

\newcommand{\rP}{\overset{\rightharpoonup }{P}}

\begin{document}
\title{A strengthening  of  a degree sequence condition for Hamiltonicity  in tough graphs}
\author{Songling Shan\footnote{Auburn University, Department of Mathematics and Statistics, Auburn, AL 36849.
		Email: {\tt szs0398@auburn.edu}.   
		Supported in part by NSF grant DMS-2345869.}
		\qquad 
		 Arthur Tanyel\footnote{Auburn University, Department of Mathematics and Statistics, Auburn, AL 36849.
		 	Email:	{\tt ant0034@auburn.edu}. }
}

\date{\today}
\maketitle

\begin{abstract}
Generalizing Chv\'atal's classic 1972 result, Ho\`ang proposed in 1995 the following conjecture, which strengthens Chv\'atal's result in terms of toughness:  Let $t\ge 1$ be a positive integer and $G$ be a $t$-tough graph on $n \ge 3$ vertices with degree   sequence  $d_1, d_2, \dots, d_n$ in non-increasing order.   Suppose  for each  $i\in [1, \lfloor\frac{n-1}{2} \rfloor]$, if
    $
    d_i \le i   \text{ and } d_{n-i+t} < n - i $ implies 
    $d_j + d_{n-j+t} \ge n$ for all $j\in [i+1,  \lfloor\frac{n-1}{2} \rfloor]$, 
    then $G$ is Hamiltonian. Ho\`ang verified the conjecture for $t=1$. In this paper, 
    we verfity the conjecture for all $t\ge 4$. Our proof relies on 
    a toughness closure lemma for $t\ge 4$ that  we previously established. Additionally, we show 
    that the toughness closure lemma does not hold when $t=1$.

\end{abstract}

\medskip

\section{Introduction}

We consider only simple graphs. Let $G$ be a graph. 
We denote by $V(G)$ and $E(G)$  the vertex set and  the edge set of $G$, respectively.
If vertices $u$ and $v$ are adjacent in $G$, we write $u \sim v$.
Otherwise, we write $u \nsim v$.
For any $S \subseteq V(G)$, let $G[S]$   be the subgraph induced on $S$, and let 
 $G-S=G[V(G) \setminus S]$. When $S=\{v\}$ is a singleton, 
 we write $G - v$ for $G - \{v\}$. 
 If $u$ and $v$ are nonadjacent vertices in $G$, we let $G + uv$ be the graph obtained by adding 
 the edge $uv$ to $G$. 
 For  $F\subseteq E(G)$, we let $G-F$ be obtained from $G$ by deleting all the edges contained in $F$. 
For integers $p$ and $q$, we write $[p,q] = \{i \in \mathbb{Z} : p \le i \le q \}$.

Let   $c(G)$ denote  the number of components of a graph $G$. 
The  \emph{toughness}  of $G$,  denoted $\tau(G)$, is $\min\{|S|/c(G-S): S\subseteq V(G), c(G-S) \ge 2\}$  if $G$ is not
a complete graph, and is defined to be $\infty$ otherwise.  A graph    is called  \emph{$t$-tough}
if  its toughness is at least $t$. 
This concept was introduced by Chv\'atal~\cite{chvatal1973tough} in 1973 as a measure of a graph's ``resilience" under the removal of vertices.

Let $n\ge 1$ be an integer. 
The non-decreasing sequence $d_1, d_2, \dots , d_n$ is a \textit{degree sequence} of graph $G$ if the vertices of $G$ can be labeled as $v_1, v_2, \dots, v_n$ such that $d(v_i) = d_i$ for all $i\in [1,n]$. In 1972, 
Chv\'atal~\cite{chvatal1972hamilton} proved the following well known result. 
\begin{THM}\label{thm:chv_degree}
	Let $G$ be a graph on $n\ge 3$ vertices  with degree sequence $d_1, d_2, \dots, d_n$. 
    If for all $i < \frac{n}{2}$, $d_i \le i$ implies $d_{n-i} \ge n-i$, then $G$ is Hamiltonian.
\end{THM}

Ho\`ang~\cite[Conjecture 1]{hoang1995hamiltonian} in 1995 conjectured a toughness analogue for the theorem above. 

\begin{CON}\label{con:degree-sequence-t-tough}
	Let  $n\ge 3$ and $t\ge 1$ be integers, and $G$ be  $t$-tough  graph with degree sequence $d_1, d_2, \dots , d_n$.  If  for all $i < \frac{n}{2}$ it holds that  $d_i \le i$ implies  $d_{n-i+t} \ge n-i$, then  $G$ is Hamiltonian.
\end{CON}

Ho\`ang in the same paper~\cite[Theorem 3]{hoang1995hamiltonian} proved the conjecture  for $t \le 3$. 
Since every hamiltonian graph must necessarily be 1-tough, the statement for $t=1$ generalizes Theorem~\ref{thm:chv_degree}. 
Recently, Ho\`ang  and Robin~\cite{hoang2024closure} proved the conjecture  for $t = 4$, and 
we  confirmed  Conjecture~\ref{con:degree-sequence-t-tough}  for all $t \ge 4$~\cite{shan2025degreesequence}.

Ho\`ang~\cite[Conjecture 4]{hoang1995hamiltonian}  also proposed  a strengthening of  Conjecture~\ref{con:degree-sequence-t-tough}.  (Ho\`ang wrote $d_j + d_{n-j+1} \ge n$ in Conjecture 4 of~\cite{{hoang1995hamiltonian}}  but  that should be a typo according to the statement of Conjecture 5 in~\cite{hoang1995hamiltonian}.)

\begin{CON}\label{con:degree-sequence-t-tough2}
Let $t\ge 1$ be a positive integer and $G$  be a $t$-tough graph on $n \ge 3$ vertices with degree sequence  $d_1, d_2, \dots, d_n$.   Suppose  for each  $i\in [1, \lfloor\frac{n-1}{2} \rfloor]$, if
    $
    d_i \le i   \text{ and } d_{n-i+t} < n - i $ implies 
    $d_j + d_{n-j+t} \ge n$ for all $j\in [i+1,  \lfloor\frac{n-1}{2} \rfloor]$, 
    then $G$ is Hamiltonian.
\end{CON}

Ho\`ang~\cite{hoang1995hamiltonian} verified the conjecture when $t=1$. We here confirm all 
the   $t\ge 4$ cases.

\begin{THM}\label{thm:degree-sequence-t-tough}
	Let $t\ge 4$ be a positive integer and $G$ be a $t$-tough graph on $n \ge 3$ vertices with degree sequence  $d_1, d_2, \dots, d_n$.   Suppose  for each  $i\in [1, \lfloor\frac{n-1}{2} \rfloor]$, if
    $
    d_i \le i   \text{ and } d_{n-i+t} < n - i $ implies 
    $d_j + d_{n-j+t} \ge n$ for all $j\in [i+1,  \lfloor\frac{n-1}{2} \rfloor]$, 
    then $G$ is Hamiltonian.
\end{THM}

The proof of Theorem~\ref{thm:degree-sequence-t-tough} relies on a toughness 
closure lemma  for toughness at least 4 that we established in~\cite{shan2025degreesequence}.  
In this paper, we also show that the toughness closure lemma does not hold for $t=1$. 
To state the result, we provide some definitions about  Hamiltonian closures. 

Let $G$ be a graph on $n$ vertices. 
 Bondy and Chv\'atal~\cite{bondychvatal1976} in 1976 defined the closure of $G$ 
as the graph obtained from $G$ by iteratively adding all edges
joining pairs of nonadjacent vertices whose degree-sum is at least $n$
in the current stage. They further proved the following classic result. 

\begin{THM}[Bondy and Chv\'atal~\cite{bondychvatal1976}]\label{thm:closure-lemma}
   A graph  $G$ is Hamiltonian if
and only if its closure is Hamiltonian. 
\end{THM}

 When investigating Conjecture~\ref{con:degree-sequence-t-tough},  Ho\`ang and Robin~\cite{hoang2024closure} in 2024 introduced 
the analogous closure concept for $t$-tough graphs $G$, called the $t$-closure, where $t\ge 1$ is 
an integer. The \emph{$t$-closure} of $G$ is 
the graph obtained from $G$ by iteratively adding all edges
joining pairs of nonadjacent vertices whose degree-sum is at least $n-t$. 
Ho\`ang and Robin~\cite{hoang2024closure}  further showed that for $t\ge 2$, a $\frac{3t-1}{2}$-tough graph $G$ is Hamiltonian if and only if its $t$-closure is Hamiltonian. 

We established  the exact toughness analogue of the Bondy-Chv\'atal  closure lemma when $t\ge 4$ 
as given below.  
In this paper,  we provide a counterexample demonstrating that this statement fails for $t=1$.

 \begin{THM}[{Shan and Tanyel~\cite[Theorem 6]{shan2025degreesequence}}]\label{thm:t-closure}
	Let $t \geq 4$ be an integer,   $G$ be  a $t$-tough graph on $n \geq 3$ vertices, and let  distinct   $x, y \in V(G)$  be nonadjacent with degree-sum at least $n-t$. 
	Then $G$ is Hamiltonian if and only if  $G+xy$ is Hamiltonian.  
\end{THM}

 \begin{THM}\label{1-closure}
 For any integer $n\ge 7$, there exists a graph $G$ on $n$ vertices with the following properties:
 \begin{enumerate}[(1)]
 
     \item There exist nonadjacent vertices $x,y \in V(G)$ such that $d(x) + d(y) = n - 1$; 
     \item $G+xy$ is Hamiltonian but $G$ is not Hamiltonian;
     \item $\tau(G)=1$. 
 \end{enumerate}
\end{THM}

The lower bound of $n$ in Theorem~\ref{1-closure} is best possible, as demonstrated below.  

 \begin{THM}\label{thm:lower-bound-n}
    For any  integer $n\in [3,6]$, if $G$
is a 1-tough $n$-vertex graph with  nonadjacent $x,y\in V(G)$
for which  $d(x)+d(y) \ge n-1$, then $G+xy$ is Hamiltonian implies 
that $G$ is Hamiltonian.
\end{THM}

The remainder of this paper is organized as follows.  In Section 2, we prove Theroem~\ref{thm:degree-sequence-t-tough} 
by applying Theorem~\ref{thm:t-closure}. In the last section, we prove Theorem~\ref{1-closure} 
and Theorem~\ref{thm:lower-bound-n}.

\section{Proof of Theorem~\ref{thm:degree-sequence-t-tough}}

Let $G$ be a graph. 
For  $S,T \subseteq V(G)$, we say that $S$ is \emph{complete to} $T$ if for all $u \in S$ and $v \in T \setminus \{u\}$, it holds that $u \sim v$.
If $S$ is complete to $V(G)$, we call $S$ a \emph{universal clique} in  $G$. We will need the two results below in our proof.

\begin{THM}[Bauer et al.~\cite{bauer1995long}]\label{bauer}
	Let  $t \ge 0$ and $G$ be a $t$-tough graph on $n \ge 3$ vertices.
    If  $\delta(G) > \frac{n}{t+1} - 1$, then  $G$ is Hamiltonian.
\end{THM}

 \begin{THM}\label{shanandtanyel1}
	Let $t \ge 4$ be an integer and  $G$ be  a $t$-tough graph on $n \geq 3$ vertices with degree sequence $d_1, d_2, \dots , d_n$.  If  for all $i < \frac{n}{2}$ it holds that  $d_i \le i$ implies  $d_{n-i+t} \ge n-i$, then  $G$ is Hamiltonian.
\end{THM}

\proof[Proof of Theorem~\ref{thm:degree-sequence-t-tough}]
Let $t\ge 4$ be an integer and $G$ a $t$-tough graph with degree sequence $d_1, d_2, \dots, d_n$ as described.  We let $v_i$ be the vertex of $G$ such that $d(v_i)=d_i$ for each $i$. 
We assume to the contrary that $G$ is not Hamiltonian.  This, in particular, implies that $G$ is not a complete graph and so $\delta(G) \ge 2t$.  
By Theorem~\ref{shanandtanyel1}, we also 
 assume that there exists an integer $h$ with 
 $$\text{$1 \le h < \frac{n}{2}$ such that $d_h \le h$ and $d_{n-h+t}< n - h$}.$$
The hypothesis of the theorem implies the following fact:
\begin{equation}\label{eqn:less-than-h-indices}
    \text{For any $i<h$ with $d_i\le i$, it holds that $d_{n-i+t} \ge n-i$}. 
\end{equation}
 
 As adding edges to $G$ preserves the same condition on the resulting degree sequence, 
 we may assume that $G$ is its $t$-closure by 
Theorem~\ref{thm:t-closure}.  Thus 
\begin{equation}\label{eqn:degree-sum}
    \text{for any two distinct vertices $u,v\in V(G)$,   $d(u)+d(v) \ge n-t$ 
implies $u\sim v$. }
\end{equation}

\begin{CLA}\label{claim:j-and-n-j-degree-bound}
  For any $i$ with $i\in[h+1,\frac{n-1}{2} \rfloor] $, it holds that $d_i \ge i-t+1$ and $d_{n-i+t} \ge n - i$.
\end{CLA}

\proof[Proof of Claim~\ref{claim:j-and-n-j-degree-bound}] 
Since $d_i+d_{n-i+t} \ge n$ by the condition in Theorem~\ref{thm:degree-sequence-t-tough}, we have
$d_i + d_j \ge n$ for any $j \ge n - i + t$ and $d_j + d_{n-i+t} \ge n$ for any $j \ge i$. 
Thus the vertex $v_i$ has at least $n - (n - i + t) + 1 = i-t+1$ neighbors
and the vertex $v_{n-i+t}$ has at least $n-1 - (i-1)=n-i$ neighbors. 
\qed

Let $k\in [1, \frac{n-1}{2} \rfloor]$ be the smallest integer such that $d_k  \le  k$. 
Then we have   $d_i > i$ for all $i\in [1,k-1]$. Thus $d_k=k$. 
Since we have $d_h\le h$ and $k$ is the smallest integer with $d_k\le k$, it follows
 that $k\le h$. 
Further, we have $k\ge 2t$ by   $\delta(G) \ge 2t$.

Our goal  below is to find a universal clique of size larger than  $\frac{n}{t+1} - 1$ in $G$.
This would imply  that $G$ is Hamiltonian by Theorem~\ref{bauer}. 

For any integer $\alpha \in [1,\frac{n-1}{2} \rfloor]$,  let
$$U^\alpha = \{v_i\in V(G) : d_i \ge n - \alpha\}.$$

\begin{CLA}\label{claim:condition-for-Uclique}
 Let  $\alpha \in [1, \frac{n-1}{2} \rfloor]$. 
 For every $i \in [1,n]$, if   $d_i\ge \alpha-t$ or  $d_i \ge i - t + 1$, then $U^\alpha $ is a universal clique in $G$. 
\end{CLA}

\proof[Proof of Claim~\ref{claim:condition-for-Uclique}] 
Assume to the contrary that $U^\alpha $ is not a universal clique. 
Then there exists $v_p\in U^\alpha$ and $v_q\in V(G)$ such that 
$v_p\not\sim v_q$. By Fact~\ref{eqn:degree-sum}, we must have $d_q <\alpha-t$.   We choose  
$q\in[1,n]$ to be maximum with the property that $v_p\not\sim v_q$. 
 By the hypothesis of this claim, we have $d_q \ge q - t + 1$.
By the maximality of $q$, we have $v_p \sim v_\ell$ for all $\ell \in [q+1,n]$. That is, $d_p \ge n - q - 1$.
However, this gives $d_p + d_q \ge n-q-1+q - t + 1 = n - t$, a contradiction to Fact~\ref{eqn:degree-sum}. 
\qed

\begin{CLA}\label{claim:U-clique-size}
For every  $\alpha \in [1, \frac{n-1}{2} \rfloor]$, we have that $d_\alpha \le \alpha  $ implies $|U^\alpha| \ge \alpha - t$.
\end{CLA}

\proof[Proof of Claim~\ref{claim:U-clique-size}] Consider first that $\alpha \ne h$. 
If $\alpha < h$, then we have  $d_{n-\alpha + t} \ge n - \alpha$ by~\eqref{eqn:less-than-h-indices}.
That is, there are at least $n -(n-\alpha + t) + 1= \alpha - t + 1$ vertices of degree $n - \alpha$.
Thus $|U^\alpha| \ge \alpha -t+1$. 
If $\alpha > h$, then we   have  $d_\alpha+ d_{n-\alpha + t} \ge n$ by the hypothesis of the theorem. 
Thus $d_{n-\alpha + t}+d_i \ge n$ for any $i\ge \alpha$. 
Thus $v_{n-\alpha + t}$ has at least $n-1-(\alpha-1)=n-\alpha$ neighbors,
and so $d_{n-\alpha + t} \ge n-\alpha$. Hence, there are at least $n -(n-\alpha + t) + 1= \alpha - t + 1$ vertices of degree $n - \alpha$, and so $|U^\alpha| \ge \alpha -t+1$.

We now consider  the case $\alpha = h$. If  $d_{h-1} \le h - 1$, then  
we have  $d_{n-(h-1)+t} \ge n - (h - 1) = n - h + 1$ by~\eqref{eqn:less-than-h-indices}. 
Thus, there are at least $n - (n-(\alpha-1)+t)  + 1= \alpha - t $ vertices of degree at least $n - \alpha$, and so $|U^\alpha| \ge \alpha-t$. 
Thus  we assume $d_{h-1} > h -1$.  This gives $d_{h-1}=h$
as $d_h\le h$. Then by Claim~\ref{claim:j-and-n-j-degree-bound}, 
we have 
 $d_{n-(h+1)+t} \ge n - h-1$.
Therefore  $v_i \sim v_{n-(h+1)+1}$ for all $i \ge h - 1$ by~\eqref{eqn:degree-sum}.
That is, $d_{n-(h+1)+1} \ge n - (h-1) + 1 = n - h + 2$.
Hence,  there are at least $n - (n-(h+1)+t) + 1 = h+2-t$ vertices of degree at least $n - h$, which gives $|U^h| \ge h - t$. 
 \qed

\begin{CLA}\label{main4}
Let $\Omega \subseteq V(G)$ be a maximum sized universal clique in $G$.  Then $|\Omega| \le k - 2$.
\end{CLA}

\proof[Proof of Claim~\ref{main4}] 
Assume $|\Omega| \ge k - 1$. 
As $\Omega$ is a universal clique, $d_i \ge |\Omega| \ge k - 1$ for all $i \in [1,n]$. 
If $|\Omega| > k $, then $d_1 > k$ which contradicts $d_1 \le d_k = k$. Thus, $|\Omega| \in [k-1,k]$. 
Observe that $v_i \notin \Omega$ whenever $i \in [1,k]$ since every element of $\Omega$ has degree $n-1$ and   $n - 1 > \frac{n}{2} > k$. 
Let
$$
S=(\bigcup_{i=1}^k N(v_i))\setminus \{v_1, \ldots, v_k\}. 
$$
Certainly $\Omega \subseteq S$. Each vertex $v_i$ for $i\in [1,k]$ has at most $k - |\Omega| \le 1$ neighbor in $\{v_1, \ldots, v_k\}\setminus \Omega$. 
Let $p$ be the number of edges in $G[\{v_1, \ldots, v_k\}]$. 
As edges in $G[\{v_1, \ldots, v_k\}]$ form a matching, we have $k\ge 2p$. 
Then,
 \begin{numcases}{|S| \le } 
 |\Omega| =k &\text{if $|\Omega| = k$}, \nonumber \\
 |\Omega| + k-2p \le 2k-2p-1  &\text{if $|\Omega| = k-1$}. \nonumber
 \end{numcases} 
 When $|\Omega| = k$ and  $p = 0$,  we have  $c(G - S) \ge k \ge 2$. 
 If $|\Omega | = k - 1$ and $p=1$, then we have $c(G-S) \ge k-p \ge 3$. 
 If $|\Omega | = k - 1$ and $p\ge 2$, then we have $c(G-S) \ge k-p \ge p\ge 2$. 
 Thus $\frac{|S|}{c(G-S)} \le \frac{2k-2p-1}{k-p}< 2<4$, contradicting the toughness of $G$. 
 \qed

  \begin{CLA}\label{main6}
For all $\alpha$ with $k + t - 1 \le \alpha < \frac{n}{2}$, we have $d_\alpha \ge \alpha-t+1$ if $\alpha >h$ and $d_\alpha > \alpha$ if $\alpha \le  h$. 
\end{CLA}

\proof[Proof of Claim~\ref{main6}] If $\alpha > h$, then we have $d_\alpha \ge \alpha-t+1$ 
by  Claim~\ref{claim:j-and-n-j-degree-bound}. 
Thus we assume  $\alpha \le h$.
Suppose to the contrary that $d_\alpha \le \alpha$.
Choose $\alpha$ to be the minimum such integer.

It suffice to show that 
$U^\alpha$ is a universal clique for  achieving  a contradiction: 
by Claims~\ref{main4} and~\ref{claim:U-clique-size}, we have $k - 2 \ge |\Omega| \ge |U^\alpha| \ge \alpha - t$. 
Rearranging the inequalities gives $k + t - 2 \ge \alpha \ge k + t - 1$, contradiction.

Thus, in the following, we show that $U^\alpha$ is a universal clique. 
By  Claim~\ref{claim:condition-for-Uclique},  we  only need to show that for all $i \in [1, n]$, either $d_i\ge \alpha-t$ or $i - t + 1 \le d_i < \alpha - t$. 

If $i \in [1,k]$, then 
by the minimality of $k$, we have $d_i \ge i \ge i - t + 1$. 
If $i \in [k, k + t - 2]$, then 
  $d_i \ge d_k = k > k - 1 \ge i - t + 1  $. 
If   $i\in [k+t-1, \alpha-1]$, then  
by the minimality of $\alpha$, we have $d_i > i > i - t + 1$. 
Consider next  that $i = \alpha$.
Then, if $\alpha > k + t - 1$, we have $\alpha - 1 \ge k + t - 1$.
By the minimality of $\alpha$, we get $\alpha - 1 < d_{\alpha - 1} \le d_\alpha \le \alpha$.
Thus $d_\alpha = \alpha \ge \alpha - t$.
If $\alpha = k + t - 1$, then $d_\alpha \ge d_k = k > \alpha - t$.
In either case, $d_i \ge  \alpha - t$.
If  $i \in [\alpha + 1,n]$, then  we have $d_i \ge d_\alpha > \alpha - t$ by the argument in the case where $i = \alpha$.
%
\qed

   \begin{CLA}\label{main7}
We have $k > \frac{n}{2} - t$.
\end{CLA}

\proof[Proof of Claim~\ref{main7}] 
Assume to the contrary that $k \le \frac{n}{2} - t$. Let $p = \lfloor \frac{n-1}{2} \rfloor$.
Then $k + t - 1 \le p < \frac{n}{2}$. As $k\le h$,  by Claim~\ref{main6}, we have $d_p \ge p$.
If $d_p = n - 1$, then $\{v_i : p \le i \le n\}$ is a universal clique in $G$ of size larger than $\frac{n}{2}$, contradicting Claim~\ref{main4}.
Thus, there exists $i \in [1,n]$ such that $v_p \nsim v_i$. Choose such $i$ to be maximum.
As $v_i \nsim v_p$, we must have $d_i < n - t - d_p < n - t - (\frac{n-1}{2} - 1) = \frac{n+1}{2} - t +1 \le d_p$, which gives $i < p$.
If $i \in [1,k]$, then 
by the minimality of $k$, we have $d_i \ge i \ge i - t + 1$. 
If $i \in [k, k + t - 2]$, then 
  $d_i \ge d_k = k > k - 1 \ge i - t + 1  $. 
If   $i\in [k+t-1, p-1]$, then  we have $d_i  \ge  i - t + 1$
by Claim~\ref{main6}. 
By the maximality of $i$, we have $v_p \sim v_j$ whenever $j \in [i+1,n]$.
Thus $d_p + d_i \ge n -1- i + i - t + 1 = n - t$, a contradiction to~\eqref{eqn:degree-sum}. 
\qed

   \begin{CLA}\label{main8}
We have $\delta(G) > \frac{n}{t+1}-1$.
\end{CLA}

\proof[Proof of Claim~\ref{main8}] 
Assume that $\delta(G) \le \frac{n}{t+1}-1$.
Then, as $2t \le \delta(G)$, we have $(2t+1)(t+1) \le n$.
By the minimality of $k$, we have $d_i>i$ for any $i\in [1,k-1]$. 
Also, we have $d_i\ge d_k =k > k-t$ for any $i\in [k,n]$. 
Thus by Claim~\ref{claim:condition-for-Uclique}, $U^k$ 
is a universal clique in $G$. 
Therefore, by Claims~\ref{claim:U-clique-size} and~\ref{main7}, we get $\delta(G) \ge |U^k| \ge k - t> \frac{n}{2} - 2t$.
Observe that for $t \ge 3$,  we have 
\begin{eqnarray*}
	\frac{n}{2}-\frac{n}{t+1} &=& \frac{n(t-1)}{2(t+1)} \ge 
 \frac{(2t+1)(t+1)(t-1)}{2(t+1)}  \\
	 &=&(t+0.5)(t-1) > 2t-1. 
\end{eqnarray*}
This gives  $ \frac{n}{2} - 2t > \frac{n}{t+1} - 1$.
Thus $\delta(G) \ge  k -t> \frac{n}{t+1} - 1$, a contradiction.
\qed

Now $G$ is Hamiltonian by  Claim~\ref{main8} and Theorem~\ref{bauer}, completing the proof. 
\qqed

\section{Proof of Theorem~\ref{1-closure}}

Let $P=v_1\ldots v_n$ be a path. 
For $u,v \in V(P)$, we let $u P v$ be the subpath of  $P$ with endvertices as $u$ and $v$. 
We use $\rP$ to denote an orientation of the path $P$.
In this paper, we assume that the orientation is in the direction of increasing indices.
We use $u^+$ to denote the immediate successor of $u$ on $P$ and $u^-$ to denote the immediate predecessor of $u$ on $P$.

	\begin{figure}[!htb]
		\begin{center}
	
			\begin{tikzpicture}[scale=1]

					{\tikzstyle{every node}=[draw ,circle,fill=black, minimum size=0.3cm, inner sep=1pt]						
						\node[draw, label=left: $x$ ] at (0,0) (x){};
						\node[draw, label=above: $\,v_2$] at (1,0) (v2){};

                          \node[draw, label=above: $\,v_{n-4}$] at (6,0) (v6){};

                           \node[draw,label=above: $v_{n-3}$] at (7,0) (v7){ };
                             \node[draw,label=above: $v_{n-2}$] at (8,0) (v8){ };
                           \node[draw,label=above: $v_{n-1}$] at (9,0) (v9){ };
                           \node[draw,label=right: $y$] at (10,0) (y){};
		
					}

                    {\tikzstyle{every node}=[draw ,circle,fill=black, minimum size=0.1cm, inner sep=1pt]						
						\node[draw ] at (1.5,0) (a){};
						\node[draw] at (2,0) (b){};
						\node[draw] at (2.5,0) (c){};

					}

					\path[draw,black]
					(x) edge node[name=la,pos=0.7, above] {\color{blue} } (v2)
                (v2) edge node[name=la,pos=0.7, above] {\color{blue} } (v6)
                    (v6) edge node[name=la,pos=0.7, above] {\color{blue} } (v9)
                    (v9) edge node[name=la,pos=0.7, above] {\color{blue} } (y)
                    ;

                    \draw [black,thick] plot [smooth,  tension=1.5] coordinates {(y)(8,-0.7) (6,0)};

                    \draw [black,thick] plot [smooth,  tension=1.5] coordinates {(y)(5.5,-1.5) (1,0)};

                    \draw [black,thick] plot [smooth,  tension=1.5] coordinates {(y)(9,-0.5) (8,0)};

                    \draw [black,thick] plot [smooth,  tension=1.5] coordinates {(x)(4, 1.6) (8,0)};

					\node[] at (5,-2) (){$N(x)=\{v_2, v_{n-2}\}$ and $N(y) =V(G)\setminus \{x,y,v_{n-3}\}$};

					\end{tikzpicture}
		
		\end{center}
        \caption{A depiction of $G$.}
        \label{f1}
	\end{figure}
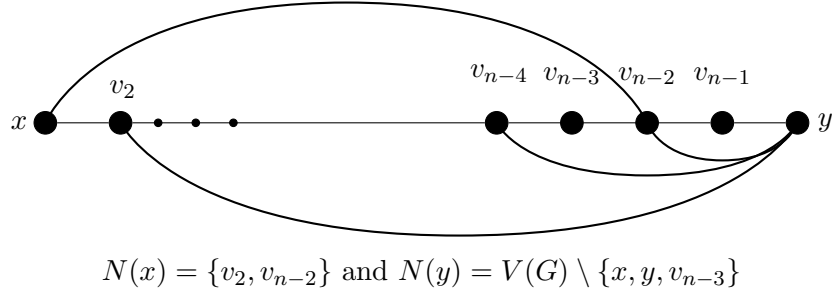

\proof[Proof of Theorem~\ref{1-closure}]
Let $n\ge 7$.  The graph $G$ has a Hamiltonian path 
$P=v_1\ldots v_n$, where $x=v_1$ and $y=v_n$, 
and $E(G)=E(P)\cup \{xv_{n-2}\} \cup \{yv_i: i\in [2,n-2]\setminus \{n-3\}\}$. See Figure~\ref{f1} for a depiction.

It is clear that $x\not\sim y$ and $d(x)+d(y)=2+n-3=n-1$, and $G+xy$
is Hamiltonian. We next show that $G$ 
is not Hamiltonian. Suppose to the contrary that $G$
has a Hamiltonian cycle $C$.  Then as $d(x)=d(v_{n-3})=d(v_{n-1})=2$, it follows 
that $xv_2, xv_{n-2}, v_{n-4}v_{n-3}, v_{n-3}v_{n-2}, v_{n-2}v_{n-1}, v_{n-1}y\in E(C)$. Howerver, 
this implies that the degree of $v_{n-2}$ in $C$ is 3, a contradiction. 

Next, we show that $\tau(G)=1$.  As $\delta(G)=2$, we have $\tau(G) \le 1$. Let $S$ be any cutset of $G$. We show that $c(G-S) \le |S|$. 
Note that $C=xv_{n-2} v_{n-1}y v_{n-4}Px$ is a  Hamiltonian cycle of $G-v_{n-3}$. 
Thus $G-v_{n-3}$ is 1-tough (any cycle is 1-tough).  

Consider first that $y\not\in S$. 
If $v_{n-3}$ is not a component of  $G-S$, then 
as $y\sim v_{n-2}$ and $y\sim v_{n-4}$,   we know that $c(G-S) =c(G-v_{n-3}-(S\setminus \{v_{n-3}\}))$. 
Thus $c(G-S) \le |S|$ by  $c(G-v_{n-3}-(S\setminus \{v_{n-3}\})) \le |(S\setminus \{v_{n-3}\}| \le |S|$.  Thus we assume that $v_{n-3}$ is a trivial component of $G-S$. This implies that $ v_{n-2}, v_{n-4}\in S$. As $Q=xPv_{n-5}yv_{n-1}$ is a path
and so $c(Q-(S\cap V(Q))) \le |S\cap V(Q)|+1$, we get that 
$|S|=|S\cap V(Q)|+2 \ge c(Q-(S\cap V(Q)))+1=c(G-S)$. 

Consider next that $y\in S$. Then $Q=P-y$ is a path, and we again 
have $c(Q-(S\cap V(Q))) \le |S\cap V(Q)|+1$. 
Thus $|S|=|S\cap V(Q)|+1 \ge c(Q-(S\cap V(Q)))=c(G-S)$. 
\qqed 

Lastly, we  prove that the lower bound of $n$ in Theorem~\ref{1-closure} is best possible.

\proof[Proof of Theorem~\ref{thm:lower-bound-n}]
Suppose to the contrary that $G$
is not Hamiltonian.  Then $G$ has a Hamiltonian path
$P=v_1\ldots v_n$, where $x=v_1$ and $y=v_n$.

The first three  assertions below follow directly from the assumption that $G$ is not Hamiltonian, and the third is a corollary of the first three. 

\begin{CLA}\label{claim:adjacency}
Let  distinct $i,j \in [2,n-1]$ and suppose $x\sim v_i$ and $y\sim v_j$. Then the following holds. 
\begin{enumerate}[(1)]
\item We have $y\not\sim v_i^-$ and $x\not\sim v_j^+$. 
	\item  If $i<j$, then $v_i^- \not\sim v_j^+$. 
	\item If $i>j$, then $v_i^+\not\sim v_j^+$ and $v_i^-\not\sim v_j^-$. 
      \item  If $v_i^-\sim y$ and $v_i^+\sim x$ for some $i\in [3,n-2]$, then 
      for any $v_j\in N(x)\cup N(y)$, we have $v_i\not\sim v^+_j$
      if $j>i$ and $v_i\not\sim v^-_j$ if $j<i$. 
\end{enumerate}
\end{CLA}
Let
$$
R=N(x)\cap N(y) \quad \text{and} \quad  S=V(G) \setminus (N(x) \cup N(y) \cup \{x,y\}). 
$$
Then we have 
\begin{equation}\label{eqn:size-of-N(x)UN(y)}
    n - 1  = d(x) + d(y)  = |N(x) \cup N(y)| + |N(x) \cap N(y)| = n-2-|S|+|R|. 
\end{equation}
The equation gives $|S| = |R|-1$. 

\begin{CLA}\label{2.4}
    $|R | \ge 2$.
\end{CLA}

\proof[Proof of Claim~\ref{2.4}] Assume to the  that $|R| \le 1$.  
Then  $|R | = 1$ 
 and   $|N(x) \cup N(y)| = n - 2$ by Equation~\eqref{eqn:size-of-N(x)UN(y)}.
That is, $N(x) \cup N(y) = V(G) \setminus \{x,y\}$. By   Claim~\ref{claim:adjacency}(1) that $y\not\sim v_i^-$ for any $v_i$
with $x\sim v_i$, it follows that $v_i^- \in N(x)$. 
Thus there  exists $k \in [2,n-1]$ such that $N(x) = \{v_1, \dots , v_k\}$ and $N(y) = \{v_{k+1}, \dots, v_{n-1}\}$. By   Claim~\ref{claim:adjacency}(2),  we know that 
$G - v_k$ has two components, contradicting that $G$ is 1-tough.
\qed

As $n\ge |R|+|S|+|\{x,y\}|$, $|S| = |R|-1$,  and $n\le 6$,  it follows that $|R| \le 2$.
Thus, $|R|=2$.  Let $u$ be the vertex in $S$. If $u^-, u^+\in R$, then 
by Claim~\ref{claim:adjacency}(2) and (4), we have  $c(G-\{u^+, u^-\})=3$, 
contradicting the toughness of $G$.  Thus $|\{u^+, u^-\}\cap R| \le 1$. 

Assume that $R=\{v_i, v_j\}$ with $i<j$ and $i,j\in [2,5]$. As $u\in V(v_iPv_j)$ 
by Claim~\ref{claim:adjacency}(1), and $|\{u^+, u^-\}\cap R| \le 1$, it follows 
that $n=6$ and  $R= \{v_2, v_{5}\}$.  By symmetry, assume $u=v_3$. Then $x\sim v_4$
by Claim~\ref{claim:adjacency}(1). Now  $C=xv_5yv_2v_3v_4 x$ is a Hamiltonian cycle of $G$. 
This gives a contradiction to our assumption that $G$ is not Hamiltonian. 
\qqed


\end{document}